\documentclass[oneside,leqno,11pt]{article}
\usepackage{amssymb,amsmath,latexsym}

\def\ga{\mathfrak{a}}

\def\gg{\mathfrak{g}}

\def\gk{\mathfrak{k}}

\def\gn{\mathfrak{n}}

\def\gp{\mathfrak{p}}
\def\gq{\mathfrak{q}}

\def\C{\mathbb{C}}

\def\P{\mathbb{P}}

\def\R{\mathbb{R}}

\def\Z{\mathbb{Z}}

\def\cB{\mathcal{B}}
\def\cC{\mathcal{C}}

\def\cO{\mathcal{O}}

\def\cS{\mathcal{S}}

\def\cU{\mathcal{U}}

\def\cX{\mathcal{X}}
\def\cY{\mathcal{Y}}

\newtheorem{theorem}[equation]{Theorem}

\newtheorem{lemma}[equation]{Lemma}

\newtheorem{corollary}[equation]{Corollary}

\newtheorem{proposition}[equation]{Proposition}

\newtheorem{definition}[equation]{Definition}

\newtheorem{remark}[equation]{Remark}

\title{Schubert Varieties and Cycle Spaces}

\author{Alan T. Huckleberry\footnote{
Research partially supported by Schwerpunkt "Global methods in complex
geometry" and SFB-237 of the Deutsche Forschungsgemeinschaft.}
\,\, \& Joseph A. Wolf\footnote{
Research partially supported by NSF Grant DMS 99-88643 and by
the SFB-237 of the Deutsche Forschungsgemeinschaft.}
}

\begin{document}

\maketitle

\abstract{\begin{quote} \footnotesize 
Let $G_0$ be a real semisimple Lie group.  It acts naturally on every
complex flag manifold $Z = G/Q$ of its complexification.  
Given an Iwasawa decomposition
$G_0 = K_0A_0N_0$\,,  a $G_0$--orbit $\gamma \subset Z$, and the dual 
$K$--orbit $\kappa \subset Z$, Schubert varieties are studied
and a theory of Schubert slices for arbitrary $G_0$--orbits
is developed. For this, certain geometric properties of dual pairs
$(\gamma ,\kappa )$ are underlined.  
Canonical complex analytic slices contained in a given $G_0$-orbit $\gamma $
which are transversal to the dual $K_0$-orbit $\gamma \cap \kappa $ are
constructed and analyzed. Associated algebraic incidence divisors
are used to study complex analytic properties of certain cycle domains.
In particular, it is shown that the linear cycle space $\Omega_W(D)$ 
is a Stein domain that contains the universally defined 
Iwasawa domain $\Omega_I$\,.  This is one of the main ingredients
in the proof that $\Omega_W(D)=\Omega _{AG}$ for all but a 
few hermitian exceptions. In the hermitian case, 
$\Omega_W(D)$ is concretely described in terms of the 
associated bounded symmetric domain.  \end{quote}
}

\setcounter{section}{-1}
\section{Introduction}   
\setcounter{equation}{0}

Let $G$ be a connected complex semisimple Lie group and $Q$ a parabolic 
subgroup. We refer to $Z = G/Q$ as a complex flag manifold.    
Write $\gg$ and $\gq$ for the respective Lie algebras of
$G$ and $Q$.  Then $Q$ is the $G$--normalizer of $\gq$. Thus we can
view $Z$ as the set of $G$--conjugates of $\gq$.  The correspondence
is $z \leftrightarrow \gq_z$ where $\gq_z$ is the Lie algebra
of the isotropy subgroup $Q_z$ of $G$ at $z$. 

Let $G_0$ be a real form of $G$, and let $\gg_0$ denote its Lie algebra.  
Thus there is a homomorphism $\varphi : G_0 \to G$ 
such that $\varphi(G_0)$ is closed in $G$ and $d\varphi : \gg_0 \to \gg$ is an
isomorphism onto a real form of $\gg$.  This gives the action of $G_0$ on $Z$.

It is well known \cite{W0} that there are only finitely many $G_0$--orbits 
on $Z$. Therefore at least one of them must be open.  

Consider a Cartan involution $\theta$ of $G_0$ and extend it as usual to $G$,
$\gg_0$ and $\gg$.  Thus the fixed point set $K_0 = G_0^\theta$ is a 
maximal compactly embedded subgroup of $G_0$ and $K = G^\theta$ is its
complexification.  This leads to Iwasawa decompositions $G_0 = K_0A_0N_0$\,.

By {\em Iwasawa--Borel subgroup} of $G$ we mean a Borel subgroup
$B \subset G$ such that $\varphi(A_0N_0) \subset B$ for some Iwasawa
$G_0 = K_0A_0N_0$\,.  Those are the Borel subgroups of the form 
$B = B_MAN$, where $N$ is the complexification of $N_0$, $A$ is the 
complexification of $A_0$, $M = Z_K(A)$ is the complexification of $M_0$\,,
and $B_M$ is a Borel subgroup of $M$.
Since any two Iwasawa decompositions of $G_0$ are
$G_0$--conjugate, and any two Borel subgroups of $M$ are 
$M_0$--conjugate because $\varphi(M_0)$ is compact, it follows that any two 
Iwasawa--Borel subgroups of $G$ are $G_0$--conjugate.  

Given an Iwasawa--Borel subgroup $B \subset G$, we study the Schubert
varieties $S = {\rm c\ell}(\cO) \subset Z$, where $\cO$ is a $B$--orbit on $Z$.
We extend the theory \cite{H} of
Schubert slices from the open $G_0$--orbits to arbitrary 
$G_0$--orbits.  As a main application we show that the corresponding Schubert
domain $\Omega_S(D)$ for an open $G_0$--orbit $D \subset Z$ is equal to
the linear cycle space $\Omega_W(D)$ considered in \cite{W1}.
That yields a direct proof that the $\Omega_W(D)$ are Stein manifolds.
Another consequence is one of the
two key containments for the complete description of linear
cycle spaces when $G_0$ is of hermitian type (Section \ref {Hermitian} and
\cite {WZ2}). The identification
$\Omega _W(D)=\Omega _S(D)$ also plays an essential role in the 
subsequently proved identification of $\Omega _W(D)$ with 
the universally defined domain $\Omega _{AG}$ \cite {FH}.

Our main technical results (Theorem \ref{schubert-slices} and 
Corollary \ref{special-slice}), which
can be viewed as being of a complex analytic nature, provide
detailed information on the $q$-convexity of $D$ and the
Cauchy-Riemann geometry of the lower-dimensional $G_0$-orbits.

We thank the referee for pointing out an error in our original argument
for the existence of supporting incidence hypersurfaces at boundary points
of the cycle space.  This is rectified in Section \ref{support} below.

\section{Duality}   
\setcounter{equation}{0}

We will need a refinement of Matsuki's $(G_0,K)$--orbit duality \cite{M}.
Write $Orb(G_0)$ for the set of $G_0$--orbits in $Z$, and similarly write 
$Orb(K)$ for the set of all $K$--orbits.  A pair 
$(\gamma, \kappa) \in Orb(G_0) \times Orb(K)$ is {\em dual} or 
{\em satisfies duality} if
$\gamma \cap \kappa$ contains an isolated $K_0$--orbit.  The duality theorem
states that 
\begin{equation}\label{mat_dual}
\text{if $\gamma \in Orb(G_0)$, or if $\kappa \in Orb(K)$, there is a
unique dual $(\gamma, \kappa) \in Orb(G_0) \times Orb(K)$.}
\end{equation}
Furthermore,
\begin{equation}\label{characterization}
\text{if $(\gamma, \kappa)$ is dual, then $\gamma \cap \kappa$ is a single 
$K_0$--orbit.}
\end{equation}
Moreover, if $(\gamma, \kappa)$ is dual, then the intersection
$\gamma \cap \kappa$ is transversal: if $z \in \gamma \cap \kappa$, then the 
real tangent spaces satisfy
\begin{equation}\label{transversal}
\text{$T_z(\gamma) + T_z(\kappa) = T_z(Z)$ and 
$T_z(\gamma \cap \kappa) = T_z(\gamma) \cap T_z(\kappa) = T_z(K_0(z))$.}
\end{equation}
We will also need a certain ``non--isolation'' property:  
\begin{align} \label{non-isolation}
\begin{tabular}{l} 
Suppose that $(\gamma, \kappa)$ is not dual but
$\gamma \cap \kappa \ne \emptyset$.  \\
If $p \in \gamma \cap \kappa$,
there exists a locally closed $K_0$--invariant submanifold
$M \subset \gamma \cap \kappa$  \\
\hskip 1in such that $p \in M$ and $\dim M = \dim K_0(p) +1$.
\end{tabular}
\end{align}


The basic duality (\ref{mat_dual}) is in \cite{M}, and the refinements 
(\ref{characterization}) and (\ref{transversal}) are given
by the moment map approach (\cite{MUV}, \cite{BL}).  See
\cite[Corollary 7.2 and \S 9]{BL}. 
\smallskip

The non--isolation property
(\ref{non-isolation}) is implicitly contained in the moment map 
considerations of \cite{MUV} and \cite{BL}. 
Following \cite{BL}, the two essential ingredients are the
following.
\begin {enumerate}
\item Endow $Z$ with a $G_u$--invariant K\" ahler metric, e.g., from the
negative of the Killing form of $\gg_u$\,.  Here
$G_u$ is the compact real form of $G$ denoted $U$ in \cite{BL}.
The $K_0$-invariant
gradient field $\nabla f^+$ of the norm function $f^+:=\Vert \mu_{K_0}\Vert ^2$
of the moment map for the $K_0$-action on $Z$, determined
by the $G_u$--invariant metric, is tangent to both the $G_0$-- and
$K$--orbits.  

\item A pair $(\gamma ,\kappa )$ satisfies duality if and only if
their intersection is non--empty and contains a point of $\{ \nabla f^+ =0\} $.
\end{enumerate}

If the pair does not satisfy duality, $p\in \gamma \cap \kappa $,
$g=g(t)$ is the $1$-parameter group associated to $\nabla f^+$  and
$\epsilon >0$ is sufficiently small, then
$$M:=\bigcup _{\vert t\vert < \epsilon }g(t)(K_0(p))$$
is the desired submanifold. For this it is only important to note
that $\nabla f^+$ does not vanish along, and is nowhere tangent to, $K_0(p)$.
That completes the argument.
\medskip

Dual pairs have a retraction property, which we prove using the
moment map approach.

\begin{theorem}\label{retract-any-orbit}
Let $(\gamma ,\kappa )$ be a dual pair.  Fix $z_0\in \gamma \cap \kappa =
K_0(z_0)$.  Then
the intersection $\gamma \cap \kappa $ is a $K_0$-equivariant
strong deformation retract of $\gamma $.
\end{theorem}

\noindent {\bf Proof.}  We use the notation of \cite[\S 9]{BL}.
If $\phi (t,x)$ is the flow of $\nabla f^+$, then
it follows that, if $x\in \gamma $ and $t>0$, then $\phi (t,x)\in \gamma $.
Furthermore, the limiting set $\pi ^+(x)= \lim_{t\to \infty } \phi (t,x)$ is
contained in the intersection $K_0(z_0)$.

Let $U$ be a $K_0$--slice neighborhood of $K_0(z_0)$ in $\gamma $.  
In other words,
if $(K_0)_{z_0}$ denotes the isotropy subgroup of $K_0$ at $z_0$,
then there is a $(K_0)_{z_0}$--invariant 
open ball $B$ in the normal space $N_{z_0}(K_0(z_0))$ such that
$U$ is the $(K_0)_{z_0}$--homogeneous fiber space  
$K_0 \times _{(K_0)_{z_0}} B$ over $B$.  The isotropy group $(K_0)_{z_0}$ 
is minimal over $U$ in that, given $z\in U$, it is $K_0$--conjugate
to a subgroup of $(K_0)_{z}$.

The flow $\phi (t,\cdot )$ is $K_0$--equivariant.  Thus, since
$\pi^+(x)\subset K_0.z_0$ for every $x\in \gamma $, every orbit
$K_0.z$ in $\gamma $ is equivariantly diffeomorphic, via some
$\phi (t_0,\cdot)$, to a $K_0$-orbit in $U$.  Consequently,
$K_0(z_0)$ is minimal in $\gamma $.

The Mostow fibration of $\gamma$ is a $K_0$--equivariant
vector bundle with total space $\gamma$ and base space that is 
a minimal $K_0$--orbit in $\gamma$.  In other words the base space is
$K_0(z_0)$.  Any such vector bundle is $K_0$--equivariantly
retractable to its $0$--section.  We may take that $0$--section to
be $K_0(z_0)$. \hfill $\square$

\begin{corollary}\label{retract-open-orbit}
Every open $G_0$-orbit in $Z$ is simply-connected.  In particular
the isotropy groups of $G_0$ on an open orbit are  connected.
\end{corollary}

Corollary \ref{retract-open-orbit} was proved by other methods in  
\cite[Theorem 5.4]{W0}. In that open orbit case of
Theorem \ref{retract-any-orbit}, $\kappa$ is the base cycle, maximal
compact subvariety of $\gamma$.  

\section{Incidence divisors associated to Schubert varieties} 
\setcounter{equation}{0}

Fix an open $G_0$--orbit $D \subset Z$.  Its dual is the
unique closed $K$--orbit $C_0$ contained in $D$.  Denote $q = \dim_\C C_0$\,.
Write $\cC^q(Z)$ for the variety of $q$--dimensional cycles in $Z$.
As a subset of $Z$, the complex group orbit $G\cdot C_0$
is Zariski open in its closure.  

At this point, for simplicity of exposition we assume that $\gg_0$ is
simple.  This entails no loss of generality because all our flags, groups,
orbits, cycles, etc. decompose as products according to the decomposition
of $\gg_0$ as a direct sum of simple ideals.

In two isolated instances of $(G_0,Z)$ (see \cite {W2}), $C_0 = Z$ and 
the orbit $G\cdot C_0$ 
consists of a single point.  If $G_0$ is of hermitian type and $D$ is an open
$G_0$--orbit of ``holomorphic type'' in the terminology of \cite{WZ1},
then $G\cdot C_0$ is the compact hermitian symmetric space dual to the
bounded symmetric domain $\cB$. This case is completely understood
(\cite{W1}, \cite{WZ1}).  In these two cases we set $\Omega:= G\cdot C_0$\,.
Here $\Omega$ is canonically identified as a coset space of $G$, because
the $G$--stabilizer of $C_0$ is its own normalizer in $G$.

Except in the two cases just mentioned, the $G$--stabilizer $\widetilde{K}$ of 
$C_0$ has identity component $K$, and there is a canonical finite equivariant 
map $\pi :G/K\to G\cdot C_0 \cong G/\widetilde{K}$.  Here we set $\Omega =G/K$.
Its base point is the coset $K$.

Suppose that $Y$ is a complex analytic subset of $Z$.
Then 
$A_Y := \{C \in \pi(\Omega) \mid \C \cap Y \ne \emptyset\}$ is a closed complex
variety in $\Omega$ \cite{BM}, called the incidence variety
associated to $Y$. 
For purposes of comparison we work with
the preimage $\pi ^{-1}(A_Y)$ in $\Omega $.  From now on we
abuse notation: we write $A_Y:=\{ C\in \Omega \mid C\cap Y\not =\emptyset \}$.  If $A_Y$ is purely of codimension $1$ then we refer to
it as the incidence divisor associated to $Y$ and denote it by
$H_Y$\,.

Now suppose that the complex analytic subset $Y$ is a Schubert variety
defined by an Iwasawa-Borel subgroup $B \subset G$.  Thus $Y$ is the
closure of one or more orbits of $B$ on $Z$.  Then the incidence variety
$A_Y$ is $B$--invariant, because $\Omega$ and $Y$ are $B$--invariant.
Define $\cY(D)$ to be the set of all Iwasawa-Schubert varieties 
$Y \subset Z$ such that
$Y \subset Z \setminus D$ \,  and $A_Y$ is a hypersurface $H_Y$\,.
Then we define the {\em Schubert domain} $\Omega_S(D)$: 
\begin{equation}\label{def_schubert_domain}
\text{$\Omega_S(D)$ is the connected component of $C_0$ in
$\Omega \setminus (\underset {Y \in \cY(D)}{\bigcup } H_Y )$.}
\end{equation}
See \cite[\S 6]{HS} and \cite{H}.  Note that any two Iwasawa-Borel
subgroups are conjugate by an element of $K_0$.  Thus

\begin {equation*}
\underset {Y\in \cY(D)}{\bigcup } H_Y=\underset {k\in K_0}{\bigcup }
k(H),
\end {equation*}
where $H:=H_1\cup \ldots \cup H_m$ is the union of the incidence
hypersurfaces defined by the Schubert varieties in the complement
of $D$ of a fixed Iwasawa-Borel subgroup.  Thus $\Omega _S(D)$
is an open subset of $\Omega $, and of course it is $G_0$-invariant
by construction.

In Corollary \ref{wolf=schubert} we will show that the cycle 
space $\Omega_W(D)$ (see (\ref{def_cyclespace})) agrees with 
$\Omega_S(D)$. Consequently, it has
the same analytic properties.  For example we now check that 
$\Omega_S(D)$ is a Stein domain.

In order to prove that $\Omega_S(D)$ is Stein, it suffices to show that
it is contained in a Stein subdomain $\widetilde{\Omega}$ of $\Omega$.  
For then, given a boundary point $p\in \text{bd}(\Omega _S(D))$ in 
$\widetilde \Omega $, it will be contained in a 
complex hypersurface $H$ that is 
equal to or a limit of incidence divisors $H_Y$\,.  Now 
$H\cap \widetilde \Omega $ is in the complement of $D$ and will be the
polar set of a meromorphic function on $\widetilde{\Omega}$.
So $\Omega_S(D)$ will be a domain of
holomorphy in the Stein subdomain  $\widetilde{\Omega}$, and 
will therefore be Stein.

As mentioned above, there are three possibilities for $\Omega$.  
If $C_0=Z$, then $\Omega$ is reduced to a point,
and $\Omega_S(D)$ is Stein in a trivial way.
Now suppose $D \subsetneqq Z$.  Then either 
$\Omega$ is a compact hermitian symmetric space $G/KP_-$ or it
is the affine variety $G/K$.
In the latter case $\Omega$ is Stein, so $\Omega_S(D)$ is Stein.
Now we are down to the case where $\Omega = G/KP_-$\, is an 
irreducible compact hermitian symmetric space.  
In particular the second Betti number $b_2(\Omega )=1$. Therefore
the divisor of every complex hypersurface in $\Omega $ is ample.
For $Y\in {\cal Y}(D)$ this implies that $\Omega \setminus H_Y$ is
affine.  Since ${\cal Y}(D)\not =\emptyset $ and 
$\Omega \setminus  H_Y\supset \Omega _S(D)$, this implies that
$\Omega _S(D)$ is Stein in this case as well.  Therefore we have
proved 
  
\begin{proposition} \label{schubert_domain_stein}
If $D$ is an open $G_0$--orbit in the complex flag manifold $Z$,
then the associated Schubert domain $\Omega_S(D)$ is Stein.
\end{proposition}
 
\section{Schubert varieties associated to dual pairs} \label{def_schu}
\setcounter{equation}{0}

Fix an Iwasawa decomposition $G_0 = K_0A_0N_0$\,.  Let $B$ be a corresponding
Iwasawa--Borel subgroup of $G$; in other words $A_0N_0 \subset B$.  
Fix a $K$--orbit $\kappa$ on $Z$ and let
$\cS_\kappa$ denote the set of all Schubert varieties $S$ defined by $B$ 
(that is, $S$ is the closure of a $B$--orbit on $Z$)
such that $\dim S + \dim \kappa = \dim Z$ and 
$S \cap \text{\,c}\ell(\kappa) \ne \emptyset$.
The Schubert varieties generate the integral homology of 
$Z$. Hence $\cS_\kappa$ is determined by the topological class of 
$\text{\,c}\ell(\kappa)$.

\begin{theorem} \label{schubert-slices} {\rm {\bf (Schubert Slices)}}
Let $(\gamma,\kappa) \in Orb(G_0) \times Orb(K)$ satisfy duality.
Then the following hold for every $S \in \cS_\kappa$\,.
\begin{itemize}
\item[{\rm 1.}]  $S \cap \text{\rm \,c}\ell(\kappa)$ is contained in 
$\gamma \cap \kappa $ and is finite.  
If $w \in S \cap \kappa $, then 
$(AN)(w) = B(w) = \cO$ where $S = \text{\rm \,c}\ell(\cO)$, and $S$ is 
transversal to $\kappa $ at
$w$ in the sense that the real tangent spaces satisfy
$T_{w}(S) \oplus T_{w}(\kappa) = T_{w}(Z)$.
\item[{\rm 2.}]  The set $\Sigma = \Sigma(\gamma,S,w) := A_0N_0(w)$ is 
open in $S$ and closed in $\gamma$.
\item[{\rm 3.}]  Let {\rm c}$\ell(\Sigma)$ and {\rm c}$\ell(\gamma)$
denote closures in $Z$.  Then the map 
$K_0 \times \text{\rm \,c}\ell(\Sigma) \to \text{\rm \,c}\ell(\gamma)$,
given by $(k,z) \mapsto k(z)$, is surjective.
\end{itemize}
\end{theorem}

\noindent {\bf Proof.} Let $w \in S \cap \text{\rm \,c}\ell(\kappa)$.
Since $\gg = \gk + \ga + \gn$, complexification of the Lie algebra
version $\gg_0 = \gk_0 + \ga_0 + \gn_0$ of $G_0 = K_0A_0N_0$\,,
we have $T_w(AN(w)) + T_w(K(w)) = T_w(Z)$.
As $w \in S = \text{\rm \,c}\ell(\cO)$ and $AN \subset B$, 
we have $\dim AN(w) \leqq \dim B(w) \leqq \dim \cO = \dim S$.
Furthermore $w \in \text{\rm \,c}\ell(\kappa)$. Thus 
$\dim K(w) \leqq \dim \kappa$.
If $w$ were not in $ \kappa$, this inequality would be strict, in violation
of the above additivity of the dimensions of the tangent spaces.
Thus $w\in \kappa $ and $T_w(S) + T_w(\kappa) = T_w(Z)$.
Since $\dim S + \dim \kappa = \dim Z$ this sum is direct, i.e.,
$T_w(S) \oplus T_w(\kappa) = T_w(Z)$.  Now also 
$\dim AN(w) = \dim S$ and 
$\dim K(w) = \dim \kappa$.  Thus $AN(w)$ is open in $S$, forcing
$AN(w) = B(w) = \cO$.  We have already seen that $K(w)$ is open 
in $\kappa $, forcing
$K(w) = \kappa $.  For assertion 1 it remains only to show that
$S \cap \kappa $ is contained in
$\gamma $ and is finite.

Denote $\widehat{\gamma} = G_0(w)$.  
If $\widehat{\gamma} \ne \gamma$, then $(\widehat{\gamma},\kappa)$ is
not dual, but $\widehat{\gamma} \cap \kappa$ is nonempty because it contains
$w$.  By the non--isolation property (\ref{non-isolation}), we have a
locally closed $K_0$--invariant manifold 
$M \subset \widehat{\gamma} \cap \kappa$ such that 
$\dim M = \dim K_0(w) +1$.  We know $T_w(S) \oplus T_w(\kappa) = T_w(Z)$,
and $K(w) = \kappa$, so $T_w(A_0N_0(w))\cap T_w(M) = 0$.  Thus
$T_w(A_0N_0(w)) + T_w(K_0(w))$ has codimension $1$ in the subspace
$T_w(A_0N_0(w)) + T_w(M)$ of $T_w(\widehat{\gamma})$, which contradicts
$G_0 = K_0A_0N_0$\,.  We have proved that $(S \cap \text{\rm \,c}\ell(\kappa))
\subset \gamma$.  Since that intersection is transversal at $w$, it is finite.
This completes the proof of assertion 1.

We have seen that $T_w(AN(w)) \oplus T_w(K(w)) = T_w(Z)$, so
$T_w(A_0N_0(w)) \oplus T_w(K_0(w)) = T_w(\gamma)$, and $G_0(w) = \gamma$.
With the characterization (\ref{characterization}) and the 
transversality conditions (\ref{transversal}) for duality we have
$\dim A_0N_0(w) = \dim T_w(\gamma) - \dim T_w(\kappa \cap \gamma)
= \dim T_w(Z) - \dim T_w(\kappa) = \dim AN(w) = \dim S$.
Now $A_0N_0(w)$ is open in $S$, . 

Every $A_0N_0$--orbit in $\gamma$ meets $K_0(w)$, because $\gamma = G_0(w)
= A_0N_0K_0(w)$.  Using (\ref{transversal}), every such $A_0N_0$--orbit
has dimension at least that of $\Sigma = A_0N_0(w)$.   Since the orbits
on the boundary of $\Sigma$ in $\gamma$ would necessarily be smaller,
it follows that $\Sigma$ is closed in $\gamma$.
This completes the proof of assertion 2.

The map $K_0 \times \Sigma \to \gamma$, by $(k,z) \mapsto k(z)$, is surjective
because $K_0A_0N_0(w) = \gamma$.  Since $K_0$ is compact and $\gamma$
is dense in $\text{\rm \,c}\ell(\gamma)$, assertion 3 follows.
\hfill $\square$
\medskip

We now apply Theorem \ref{schubert-slices} to construct an Iwasawa-Schubert
variety $Y$ of codimension $q+1$\,, $q = \dim C_0$\,, which contains
a given point $p \in \text{ bd}(D)$ and which is contained in $Z \setminus D$.
Due to the presence of the large family of $q$--dimensional cycles in
$D$, one could not hope to construct larger varieties with these
properties.

Before going into the construction, let us introduce some convenient
notation and mention several preliminary facts.

We say that a point $p \in \text{bd}(D)$ is {\em generic}, written
$p\in \text{ bd}(D)_{gen}$, if $\gamma_p := G_0(p)$ is open
in $\text{bd}(D)$.  This is equivalent to $\gamma_p$ being
an isolated orbit in $\text{bd}(D)$ in the sense that no other
$G_0$-orbit in $\text{bd}(D)$ has $\gamma_p$ in its closure.
Clearly $\text{bd}(D)_{gen}$ is open and dense in $\text{bd}(D)$.

Given $p\in \text{bd}(D)_{gen}$ the orbit $\gamma =\gamma_p$
need not be a real hypersurface in $Z$.  For example, 
$G_0 = SL_{n+1}(\mathbb R)$ has exactly two orbits in
$\mathbb P_n(\mathbb C)$, an open orbit and its complement
$\mathbb P_n(\mathbb R)$.  Nevertheless, for any $z$ in such
an orbit $\gamma $ it follows that
$c\ell (D)\cap \text {bd}(D)=\gamma $
near $z$. If $\kappa $ is dual to $\gamma $, then, since the intersection
$\kappa \cap \gamma $ is transversal in $Z$, it follows that
$\kappa \cap D\not=\emptyset $.

We summarize this as follows.
 
\begin {lemma} \label {kappa_dimension}
For $p \in \text{\rm bd}(D)_{gen}$, $\gamma =\gamma_p = G_0(p)$ and
$\kappa$ dual to $\gamma$, it follows that $\kappa \cap D \not=\emptyset$.
Furthermore, if $C_0$ is the base cycle in $D$, then
\begin {equation*}
q = \dim\ C_0 < \dim \ \kappa .
\end{equation*}
\end{lemma}

\noindent {\bf Proof.}
The property $\kappa \cap D\not=\emptyset $ has been verified above.
For the dimension estimate note that $C_0$ is dimension-theoretically
a minimal $K_0$--orbit in $D$, e.g., the $K_0$--orbits in $\kappa \cap D$
are at least of its dimension.  Since $\kappa $ is not compact,
it follows that $\dim \ \kappa > \dim \ C_0$\,.
\hfill $\square$
\medskip
 
We will also make use of the following basic fact about
Schubert varieties.
 
\begin{lemma} \label {increasing_size_of_S}
Let $B$ be a Borel subgroup of $G$, let
$S$ be a $k$-dimensional $B$--Schubert variety in $Z$, and suppose that
$\dim Z \geqq \ell \geqq k$.  Then there exists
a $B$-Schubert variety $S'$ with $\dim S' = \ell$ and $S'\supset S$.
\end{lemma}
 
\noindent {\bf Proof.}
We may assume that $S \neq Z$.  Let $\cO$ be the open $B$--orbit in $S$
and $\cO'$ be a $B$--orbit of minimal dimension among those
orbits with c$\ell (\cO') \supsetneqq \cO$.
 
For $p\in \cO$ it follows that c$\ell (\cO')\setminus \cO = \cO'$ near
$p$. Since $\cO'$ is affine, it then follows that
$\dim \ \cO' = (\dim \ \cO)+1$.  Applying this argument recursively, we find
Schubert varieties $S':= \text{c}\ell (\cO')$ of every intermediate dimension
$\ell $.
\hfill $\square$
\medskip
 
We now come to our main application of 
Theorem \ref{schubert-slices}.

\begin{corollary}\label{special-slice}
Let $D$ be an open $G_0$--orbit on $Z$ and fix a boundary point
$p \in \text{\rm bd\,}D$.  Then there exist an Iwasawa decomposition
$G_0 = K_0A_0N_0$\,, an Iwasawa--Borel subgroup $B \supset A_0N_0$\,,
and a $B$--Schubert variety $Y$, such that
{\rm (1)} $p \in Y \subset Z \setminus D$, {\rm (2)}
{\rm codim}$_Z Y = q+1$, and {\rm (3)}
$A_Y$ is a $B$--invariant analytic subvariety of $\Omega$.
\end{corollary}

\noindent {\bf Proof.}
Let $p\in \text{bd}(D)_{gen}$, let $\gamma =\gamma_p$\,, and let $\kappa$
be dual to $\gamma$.  First consider the case where
$p\in \gamma \cap \kappa $.  From Lemma \ref{kappa_dimension},
codim$S \geqq q+1$ for every $S\in {\cal S}_\kappa $.
 
Now, given $S\in {\cal S}_\kappa $, further specify $p$ to be in
$S\cap \kappa $ and let $Y$ be a $(q+1)$-codimensional
Schubert variety containing $S$ (see Lemma \ref{increasing_size_of_S}).

Since ${\rm dim}_\C C_0=q$ and ${\rm codim}_\C Y = q+1$, if 
$Y \cap D \ne \emptyset $, then there would be a point of intersection 
$z \in Y \cap C_0$\,. Since $A_0N_0(z)\subset Y$, 
it would follow that 
\begin {equation*}
q = \dim_\C C_0 = \text{ codim}_\C\ A_0N_0(z) \geqq 
\text{ codim}_\C\ Y = q+1.
\end {equation*}
Thus $Y$ does not meet $D$.  On the
other hand, using Theorem \ref{schubert-slices}, it meets every $G_0$-orbit 
in c$\ell(\gamma )$.
Thus, by conjugating appropriately, we have the desired result for any point
in the closure of $\gamma $.  Since $\gamma $ was chosen to be an
arbitrary isolated orbit in bd$(D)$, the result follows for every
point of bd$(D)$. \hfill $\square $

\section{Supporting hypersurfaces at the cycle space boundary} \label{support}
\setcounter{equation}{0}
\medskip

Let $D  = G_0(z_0)$ be an open $G_0$--orbit on $Z$.  Let $C_0$ denote
the base cycle $K_0(z_0) = K(z_0)$ in $D$, i.e., the
dual $K$-orbit $\kappa$ to the open $G_0$-orbit $\gamma =D$.
Then the cycle space of $D$ is given by
\begin {equation} \label{def_cyclespace}
\Omega _W(D):= \text{ component of } C_0 \text{ in } 
\{gC_0 \mid g \in G \text{ and } gC_0 \subset D\}.
\end {equation}
Since $D$ is open  and $C_0$ is compact, the cycle space $\Omega_W(D)$ 
initially sits as an open submanifold of the complex 
homogeneous space $G/\widetilde{K}$, where $\widetilde K$ is
the isotropy subgroup of $G$ at $C_0$\,.  In the Appendix, 
Section \ref{lifting}, we will see (with the few Hermitian exceptions which 
have already been mentioned) that the finite covering 
$\pi : \Omega := G/K \to \widetilde \Omega :=G/\widetilde{K}$
restricts to an equivariant biholomorphic diffeomorphism of
the lifted cycle space
\begin {equation} \label{def_liftedcyclespace}
\widetilde{\Omega _W(D)}:= \text{ component of } gK \text{ in }
\{gK \mid g\in G \text{ and } gC_0 \subset D\} \subset \Omega 
\end {equation}
onto $\Omega _W(D)$.
\medskip

The main goal of the present section is, given 
$C \in \text{bd}(\Omega _W(D))$, to determine a particular
point $p\in C$ which is contained in a Iwasawa--Schubert variety
$Y$ with codim$_ZY=q+1$, so that $Y\cap D=\emptyset$,
and $A_Y = H_Y$\, is of pure codimension $1$. It is then an 
immediate consequence
(see Corollary \ref {wolf=schubert}) that $\Omega _W(D)=\Omega _S(D)$.
\medskip

Given $p\in C \cap \text{bd}(D)$, we consider Iwasawa--Schubert varieties
$S = {\rm c}\ell(\cO)$ of minimal possible dimension
that satisfy the following conditions:

\begin {enumerate}
\item $p\in S\setminus \cO := E$
\item $S\cap D\not=\emptyset $
\item The union of the irreducible components of $E$ that contain
$p$ is itself contained in $Z\setminus D$.  
\end {enumerate}

Notation: Let $A$ denote 
the union of all the irreducible components of $E$ contained
in $Z\setminus D$ and let $B$ denote the union of the remaining components
of $E$.  In particular $E=A\cup B$.
\medskip

Note that by starting with the Schubert variety $S_0:=Y$ as in
the proof of Corollary \ref{special-slice},  and by considering a chain 
$S_0\subset S_1\subset \ldots $
with $\dim\,S_{i+1}=\dim\,S_i +1$\,, we eventually come to a Schubert
variety $S=S_k$ with these properties.  Of course, given $p$, the
Schubert variety $S$ may not be unique, but $\dim\,S=:n-q+\delta \geqq n-q$.
\medskip

The following Proposition gives a constructive method for determining
an Iwasawa--Borel invariant incidence hypersurface
that contains $C$ and is itself contained in the complement
$\Omega \setminus \Omega _W(D)$.  Here $S$ is constructed as above.

\begin {proposition} \label {lowering dimension}
If $\delta >0$, then $C\cap A\cap B\not=\emptyset $.
\end {proposition}

Given Proposition \ref{lowering dimension}, take a point 
$p_1\in C\cap A\cap B$ and 
replace $S$ by a component $S_1$ of $B$ that contains $p_1$\,.
Possibly there are components of $E_1:=S_1\setminus \cO_1$ that
contain $p_1$ and also have non-empty intersection with $D$.
If that is the case, we replace $S_1$ by any such component.
Since this $S_1$ still has non-empty intersection with the
Iwasawa--Borel invariant $A$, at least some of the components
of its $E_1$ do not intersect in this way.  Continuing in
this way, we eventually determine an $S_1$ that satisfies all
of the above conditions at $p_1$.  The procedure stops because
Schubert varieties of dimension less than $n-q$ have empty intersection
with $D$.

\begin {corollary}
If $S_0$ satisfies the above conditions at $p_0$, then there
exist $p_1\in \text{\rm bd}(D)$ and a Schubert subvariety
$S_1\subset S_0$ that satisfies these conditions at $p_1$ and
has dimension $n-q$.
\end {corollary}

\noindent {\bf Proof.}
We recursively apply the procedure indicated above until
$\delta =0$. \hfill $\square$
\medskip

\begin {corollary}
If $C\in \text{\rm bd}(\Omega _W(D))$, there exists an
Iwasawa--Schubert variety $S$ of dimension $n-q$ such that
$E:=S\setminus \cO$ has non-empty intersection with $C$.
\end {corollary}

\begin{corollary}
Let $C\in \text{\rm bd}(\Omega _W(D))$.  Then there exists an 
Iwasawa--Borel subgroup $B \subset G$, and a component $H_E$ of a
$B$--invariant incidence variety $A_E$, where $E = S \setminus \cO$ as above,
such that $H_E \subset \Omega \setminus \Omega _W(D)$, $C\in H_E$\,, and
$\text{\rm codim}_\Omega H_E=1$, i.e., $H_E$ is an incidence divisor.
\end {corollary}

\noindent {\bf Proof.}
The hypersurface $E$ in $S$ is the support of an ample divisor
\cite {HS}.  Thus the trace-transform method
(\cite {BK}, see also \cite[Appendix] {HS}) produces a 
meromorphic function on $\Omega $ with a pole at $C$ and polar
set contained in $A_E$.  Hence $A_E$ has a component $H_E$ as
required. \hfill $\square$
\medskip

In the language of Section \ref{def_schu} this shows that for every
$C\in \text{bd}(\Omega _W(D))$ there exists $Y\in \cY(D)$
such that $C\in H_Y$\,.  In other words, every such boundary point is
contained in the complement of the Schubert domain $\Omega _S(D)$.
By definition $\Omega _W(D)\subset \Omega _S(D)$\,.  Using that,
the equality of these domains follows immediately:

\begin {corollary} \label{wolf=schubert}
$\Omega _W(D)=\Omega _S(D).$
\end {corollary}  

Let us now turn to certain technical preparations for the
proof of Proposition \ref {lowering dimension}. For $S$ as in 
its statement, 
let $\cU_S$ be its preimage in the universal family 
$\cU$ parameterized by $\Omega $.  The mapping 
$\pi : \cU_S\to \Omega $ is proper and surjective and the
fiber $\pi ^{-1}(C)$ over a point $C\in \Omega $ can be identified
with $C\cap S$.  All orbits of the Iwasawa--Borel group that
defines $S$ are transversal to the base cycle $C_0$\,; in
particular, $C_0\cap S$ is pure--dimensional with
$\dim\, C_0\cap S=\delta $.  Thus the generic cycle in $\Omega $
has this property.
\medskip

Choose a $1$-dimensional (local) disk $\Delta $ in $\Omega $
with $C$ corresponding to its origin, such that
$I_z := \pi ^{-1}(z)$ is $\delta $--dimensional for $z\not= 0$.
Define $\cX$ to be the closure of $\pi ^{-1}(\Delta \setminus
\{0\})$ in $\cU_S$.  The map $\pi _\cX:=\pi \vert_\cX:\cX \to \Delta $  
is proper and its fibers are purely $\delta $--dimensional.  

In the sequel we use the standard moving lemma of intersection theory
and argue using a desingularization $\widetilde \pi:\widetilde S\to S$,
where only points $E$ are blown up. Let $\widetilde E$, $\widetilde A$ and
$\widetilde B$ denote the corresponding $\widetilde \pi $-preimages.  By taking
$\Delta $ in generic position we may assume that for $z\not =0$
no component of $I_z$ is contained in $E$.  Hence we may lift the
family $\cX\to \Delta $ to a family $\widetilde{\cX}\to \Delta $
of $\delta $-dimensional varieties such that $\widetilde{\cX}\to {\cal X}$
is finite to one outside of the fiber over $0\in \Delta $.  Let
$\widetilde I_z$ denote the fiber of $\widetilde {\cal X}\to \Delta $ 
at $z \in \Delta$, and
shrink $\widetilde {\cal X}$ so that $\widetilde I:=\widetilde I_0$ is 
connected.  Since $\widetilde I_z\cap \widetilde A=\emptyset$ for
$z\not =0$, it follows that the intersection class 
$\widetilde I.\widetilde A$ in the homology of $\widetilde S$ is zero.
\smallskip

An irreducible component of $\widetilde{I}$ is one of the following types:
it intersects $\widetilde A$ but not $\widetilde B$, 
or it intersects both $\widetilde A$ and $\widetilde B$, or
it intersects $\widetilde B$ but not $\widetilde A$. Write $\widetilde{I} = 
\widetilde{I_A} \cup \widetilde{I_{AB}}\cup \widetilde{I_B}$ correspondingly.

\begin{lemma} \label{moving}
$\widetilde{I_{AB}} \not= \emptyset$.
\end{lemma}

\noindent {\bf Proof.}
Since $\widetilde {I_A}\cup \widetilde {I_{AB}}\not= \emptyset $, 
it is enough to consider
the case where $\widetilde {I_A}\not =\emptyset $.  Let $H$ be a hyperplane
section in $Z$ with $H\cap S=E$ (see e.g. \cite {HS}) and put
$H$ in a continuous family $H_t$ of hyperplanes with $H_0 = H$ such
that $H_t\cap I_A$ is $(\delta -1)$--dimensional for $t\not =0$ and
such that the lift $\widetilde E_t$ of $E_t:=H_t\cap S$ contains
no irreducible component of $\widetilde {I_A}$.  In particular,
$\widetilde E_t.\widetilde {I_A}\not =0$ for $t\not =0$.  Since
$\widetilde {I_A}.\widetilde A=\widetilde {I_A}.\widetilde E=
\widetilde {I_A}.\widetilde E_t$, it follows that
$\widetilde {I_A}.\widetilde A\not =0$.  But
$0=\widetilde I.\widetilde A=\widetilde {I_A}.\widetilde A
+\widetilde{I_{AB}}.\widetilde A$ and therefore
$\widetilde {I_{AB}}\not =\emptyset $. 
$\hfill \square$
\medskip

\noindent
{\bf Proof of Proposition \ref {lowering dimension}.} We first consider
the case where $\delta \geqq 2$. Since $\widetilde {I_{AB}}\not =
\emptyset $, it follows that some irreducible component $I'$ of $I$ 
has non-empty intersection with both $A$ and $B$.  
Of course $I'\cap E =(I'\cap A)\cup (I'\cap B)$.
But $E$ is the support of a hyperplane section, and since $\dim\,I'\geqq 2$,
it follows that $(I'\cap E)$ is connected.  In particular 
$(I'\cap A)$ meets $(I'\cap B)$.  Therefore $I'\cap A\cap B\not= \emptyset$
and consequently $C\cap A\cap B\not=\emptyset$.

\medskip
Now suppose that $\delta =1$, i.e., that $\widetilde I$ is
$1$-dimensional. 
Since $\widetilde I.\widetilde A=0$, the (non-empty) intersection
$\widetilde I\cap \widetilde A$ is not discrete.
We will show that some component of $\widetilde {I_{AB}}$ is 
contained in $\widetilde A$.  It will follow immediately that
$C\cap A\cap B\not =\emptyset $.  For this we assume to the contrary
that every component of $\widetilde I$ which is contained in 
$\widetilde A$ is in $\widetilde {I_A}$\,.  We decompose 
$\widetilde I=\widetilde I_1\cup \widetilde I_2$\,, where $\widetilde I_1$
consists of those components of $\widetilde I$ which are contained
in $\widetilde A$ and $\widetilde I_2$ of those which have discrete
or empty intersection with $\widetilde A$.

\smallskip
Now $\widetilde I_1.\widetilde A=\widetilde I_1.\widetilde E$. 
Choosing $H_t$ as above, we have 
$\widetilde I_1.\widetilde E=\widetilde I_1.\widetilde E_t\geqq 0$ for
$t\not =0$. If $\widetilde I_2\not =\widetilde I_B$, then
$\widetilde I_2.\widetilde A>0$.  This would contradict
$0=\widetilde I.\widetilde A=\widetilde I_1\widetilde A+
\widetilde I_2.\widetilde A$.  
Thus $\widetilde I_2=\widetilde I_B$ and 
$\widetilde I_1=\widetilde I_A$.  But $\widetilde I_A$ and
$\widetilde I_B$ are disjoint, contrary to $\widetilde I$ being
connected.  Thus it follows that $\widetilde {I_{AB}}$ does
indeed contain a component that is contained in $\widetilde A$.
The proof is complete.$\hfill \square $ 

\begin{remark} {\em In the non--hermitian case, the main result of \cite {FH} 
leads to a non--constructive, but very short, proof of the
existence of an incidence hypersurface $H\subset G/\widetilde K$ containing
a given boundary point $C\in \text{bd}(\Omega _W(D))$.
(The analogous construction in the Hermitian case is somewhat easier;
see Section \ref {Hermitian} below.) 
For this, note that if $S$ is a $q$--codimensional Schubert variety
with $S\cap C_0\not=\emptyset $, then, using \cite {BK} as
above, for $Y:=S\setminus \cO$, it follows that $H:=H_Y$ is indeed
a hypersurface.  Now let $\Omega _H$ be the connected component
containing the base point of 
$\Omega \setminus \bigcup_{k\in K} k(H)$.
It is shown in \cite {FH} that $\Omega _H$ agrees with the Iwasawa
domain $\Omega _I$
which can be defined as the intersection of all $\Omega _H$, where
$H$ is a hypersurface in $\Omega $ which is invariant under some
Iwasawa--Borel subgroup of $G$ (see Section \ref {Iwasawa}). 
It follows that $\Omega _S(D)=\Omega _H$,
because by definition $\Omega _I\subset \Omega _S(D)\subset \Omega _I$. 
\smallskip

By Corollary \ref{special-slice} there is an incidence variety $A_Y$ in 
$Z \setminus D$ which contains $C$ and is invariant by some Iwasawa--Borel
subgroup $B$.  Since the open $B$--orbit in $\Omega $ is affine,
there exists a $B$--invariant hypersurface $H$ which contains $A_Y$.
Since $\Omega _S(D)=\Omega _I$, it is also contained in $Z/D$ and
thus it has the desired properties. That completes the short proof.
In fact it follows that $\Omega _W(D)=\Omega _I=\Omega _S(D)=\Omega _H$\,.
\hfill $\square $
}
\end{remark}

\section{Intersection properties of Schubert slices}
\setcounter{equation}{0}

Let $(\gamma ,\kappa )$ be a dual pair and $z_0 \in \gamma \cap \kappa$. 
Let $\Sigma $ be the Schubert slice at $z_0$\,, 
i.e., $\Sigma = A_0N_0(z_0)\subset \gamma $.  
In particular
$z_0 \in \Sigma \cap \kappa$.  We take a close look at the
intersection set $\Sigma \cap \kappa$.
\medskip

Let $L_0$ denote the isotropy subgroup $(G_0)_{z_0}$,
and therefore $(K_0)_{z_0} = K_0 \cap L_0$ and
$(A_0N_0)_{z_0} = (A_0N_0) \cap L_0$\,.  Define 
$ \alpha : (K_0)_{z_0} \times (A_0N_0)_{z_0} \to L_0 $
by group multiplication.

\begin{lemma}\label{iso-prod-structure}
The map $\alpha : (K_0)_{z_0} \times (A_0N_0)_{z_0} \to L_0$
is a diffeomorphism onto an open subgroup of $L_0$.
\end {lemma}

\noindent {\bf Proof.} Since
$\dim\, K_0(z_0) + \dim\,(A_0N_0)(z_0) = \dim\,G_0(z_0)$
and $\dim\, K_0 + \dim\,(A_0N_0) = \dim\,G_0$
we have
$\dim\, (K_0)_{z_0} + \dim\, (A_0N_0)_{z_0} = \dim\, L_0$\,.
Thus the orbit of the neutral point under the action of the 
compact group $(K_0)_{z_0}$ is the union of certain components of
$L_0/(A_0N_0)_{z_0}$, i.e., 
$\text{Image\,}(\alpha )=(K_0)_{z_0} \cdot (A_0N_0)_{z_0}$ is 
an open subgroup of $L_0$\,.

The injectivity of $\alpha $ follows from the fact that
$G_0$ is the topological product $G_0 = K_0 \times (A_0N_0)$.  This product
structure also yields the fact that $(K_0)_{z_0}(m)$ is transversal to
$(A_0N_0)_{z_0}$ at every $m \in (A_0N_0)_{z_0}$\,.  Thus, $\alpha $ is a 
local diffeomorphism along $(A_0N_0)_{z_0}$ and by equivariance is therefore
a diffeomorphism onto its image. \hfill $\square$

\begin {corollary}\label{one-point-intersection}
If $L_0$ is connected, in particular if $\gamma $ is simply connected,
then $\Sigma \cap \kappa =\{ z_0\}$.
\end {corollary}

\noindent {\bf Proof.} If $z_1\in \Sigma \cap \kappa$, then
there exists $k\in K_0$ and $an\in A_0N_0$ so that 
$k^{-1}(z_0) = (an)(z_0)$, i.e., $kan\in L_0$\,.  Therefore $k\in (K_0)_{z_0}$,
$an\in (A_0N_0)_{z_0}$, and $z_1 = z_0$. \hfill $\square$

\begin {theorem}\label{open_case}
Let $D$ be an open $G_0$-orbit in $Z$, $C_0 \subset D$ the base cycle, 
$z_0 \in C_0$\,, and $\Sigma = A_0N_0(z_0)$ a Schubert slice at $z_0$\,.
If $C \in \Omega_W(D)$, then $\Sigma \cap C$ consists of a single point,
and the intersection $\Sigma \cap C$ at that point is transversal.
\end {theorem}

{\bf Proof.} 
Let $S := \text{c}\ell B(z_0)$ be the Schubert variety containing $\Sigma$,
and let $k$ denote the intersection number $[S]\cdot [C_0]$.  We know from
Theorem \ref{schubert-slices} that intersection points occur only in the
open $A_0N_0$--orbits in $S$.  
The open $G_0$--orbit $D$ is simply connected, and therefore 
Corollary \ref{one-point-intersection} applies.  Thus $\Sigma \cap C_0
= \{z_0\}$.  From Theorem \ref{schubert-slices}(1) it follows that intersection
this transversal. Hence it contributes exactly $1$ to $[S]\cdot [C_0]$.  Now we
have $k$ different open $A_0N_0$--orbits in $S$, each of which 
contains exactly one (transversal) intersection point.

Cycles $C \in \Omega_W(D)$ are homotopic to $C_0$\,. Thus $[S]\cdot [C] = k$.
As we homotopy $C_0$ to $C$ staying in $D$, 
the intersection points of course move around,
but each stays in its original open $A_0N_0$--orbits in $S$.  Since
$\Sigma$ is one of those open $A_0N_0$--orbits, it follows that 
$\Sigma \cap C$ consists of a single point, and the intersection there
is transversal, as asserted.  \hfill $\square$
\medskip

\begin{remark} {\em One might hope that the orbit $D$ would be
equivariantly identifiable with a bundle of type
$K\times _{(K_0)_{z_0}} \Sigma$,
but the following example shows that this is not the case.
Let $Z = \P_2(\C)$ be equipped with the standard
$SU(2,1)$-action.  Let $D$ be the open $SU(2,1)$--orbit consisting of
positive lines, i.e, the complement of the closure of the unit ball
$B$ in its usual embedding.  The Schubert slice $\Sigma$ for $D$ is
contained in a projective line tangent to bd$(B)$; see \cite{HW}.
If $z_0\in C_0\subset D$, the only $(K_0)_{z_0}$--invariant line
in $\P_2(\C)$ that contains $z_0$ and is not contained in $C_0$
is the line determined by $z_0$ and the $K_0$-fixed point in $B$.}
\end {remark}

\section{The domains $\Omega_I$ and $\Omega_{AG}$} \label {Iwasawa}
\setcounter{equation}{0}

The Schubert domain $\Omega_S(D)$ is defined as a certain subspace
of the cycle space $\Omega$.  When $G_0$ is of hermitian
type and $\Omega$ is the associated compact hermitian symmetric space,
the situation is completely understood \cite{W1}: 
$\Omega_W(D)$ is the bounded symmetric domain dual to $\Omega$ in the
sense of symmetric spaces.  Now we put that case aside.  Then 
$\Omega \cong G/K$, and we have 
\begin{equation}\label{colocation}
\Omega_W(D) = \Omega_S(D) \subset \Omega = G/K .
\end{equation}

Let $B$ be an Iwasawa--Borel subgroup of $G$.  It has only finitely
many orbits on $\Omega$, and those orbits are complex manifolds.  The
orbit $B(1K)$ is open, because $ANK$ is open in $G$, and its complement
$S \subset \Omega$ is a finite union $\bigcup H_i$ of $B$--invariant
irreducible complex hypersurfaces.  For any given open $G_0$--orbit, 
some of these $H_i$ occur in the definition (\ref{def_schubert_domain}) 
of $\Omega_S(D)$.  The {\em Iwasawa domain} $\Omega_I$ is defined as in
(\ref{def_schubert_domain}) except that we use all the $H_i$\,: 
\begin{equation}\label{def_iwasawa_domain}
\Omega_I \text{ is the connected component of }
	C_0 \text{ in } = \Omega \setminus  
	(\underset{g\in G_0}{\bigcup} g(S) ).
\end{equation}
This definition is independent of choice of $B$ because any two
Iwasawa--Borel subgroups of $G$ are $G_0$--conjugate. Just as
in the case of the Schubert domains, we note here that
\begin {equation*}
\underset{g\in G_0}{\bigcup }g(S)=\underset{k\in K_0}{\bigcup }k(S)
\end {equation*}
is closed.  By definition,
$\Omega_I \subset \Omega_S(D)$ for every open $G_0$--orbit $D$ in $Z$.

The argument for $\Omega_S(D)$ also shows that $\Omega_I$ is a Stein domain
in $\Omega$.  See \cite{H} for further properties of $\Omega_I$

The Iwasawa domain has been studied by several authors from a completely
different viewpoint and with completely different definitions.  See
\cite{B}, \cite{BLZ}, \cite{BR}, \cite{GM} and \cite{KS}.  Here is the 
definition in \cite{B}.  Let $X_0$ be the closed $G_0$--orbit in $G/B$ and
let $\cO_{max}$ be the open $K$--orbit there.  The {\em polar} 
$\widehat {X_0}$ of $X_0$ is
the connected component of $1.K$ in 
$\{gK \in \Omega \mid g \in G \text{ and } g^{-1}X_0 \subset \cO_{max}\}$.

\begin {proposition} \label{polar-iwasawa} {\em \cite {Z}}
$\widehat{X_0} = \Omega_I$\,.
\end{proposition}

\noindent {\bf Proof.} Let $\pi: G \to G/K = \Omega$ denote the projection.
As $S$ is the complement of $B\cdot K$ in $\Omega$, 
$\pi^{-1}(\Omega_I)$ is the interior of $I:= {\bigcap}_{g \in G_0} g(ANK)$.
Note that $h \in I \Leftrightarrow g^{-1}h \in ANK$ for all $g \in G_0$
$\Leftrightarrow h^{-1}g \in KAN$ for all $g \in G_0$
$\Leftrightarrow h^{-1}G_0 \subset KAN$.  Viewing $1B$ as the base point
in $X_0$\,, the condition for $hK$ being in $\widehat {X_0}$ is that
$h^{-1}G_0B \subset KB = KAN$.  Thus $h \in I \Leftrightarrow 
hK \in \widehat{X_0}$\,, in other words $\widehat{X_0} = \pi(I)
= \Omega \setminus \bigcup_{g \in G_0} g(S)$. \hfill $\square$

\begin{corollary}
The polar $\widehat{X_0}$ to the closed $G_0$--orbit $X_0$ is a Stein 
subdomain of $\Omega$.
\end{corollary}

Now we turn to the domain $\Omega_{AG}$\,.
The Cartan involution $\theta$ of $\gg_0$ defines the usual Cartan
decomposition $\gg_0 = \gk_0 + \gp_0$ and the compact real form
$\gg_u = \gk_0 + \sqrt{-1}\ \gp_0$ of $\gg$.  Let $G_u$ be the corresponding
compact real form of $G$, real--analytic subgroup for $\gg_u$\,, acting on
$\Omega = G/K$.  Then
$$
\Omega _{AG} := 
\{x \in \Omega \mid \text{ the isotropy subgroup } (G_0)_x
\text{ is compactly embedded}\}^0,
$$
the topological component of $x_0 = 1K$.
It is important to note that the action of $G_0$ on the
Akhiezer-Gindikin domain $\Omega_{AG}$ is proper
\cite{AG}.

In work related to automorphic forms (\cite{BR}, \cite{KS}) it was shown that
$\Omega_{AG} \subset \widehat{X_0}$, when $G_0$ is a classical group.  Other
related results were proved in \cite{GM}.

By means of an identification of $\widehat{X_0}$ with 
a certain maximal domain $\Omega_{adpt}$ for the adapted complex structure 
inside the real tangent bundle of $G_0/K_0$, and using basic properties of
plurisubharmonic functions, it was shown by the first author that 
$\Omega_{AG} \subset \Omega _I$ in complete generality \cite{H}.
Barchini proved the opposite inclusion in \cite{B}.  Thus
$\Omega_{AG} = \Omega _I$\,.  
In view of Theorem \ref{polar-iwasawa}, we now have

\begin{theorem}\label{i-ag}
$\Omega_I = \widehat{X_0} = \Omega_{AG}$\,.
\end{theorem}

\begin{remark}\label{rem1} 
{\em In particular, this gives yet another proof that
$\Omega_{AG}$ is Stein.  That result was first proved in \cite{BHH} where
a plurisubharmonic exhaustion function was constructed.} 
\end{remark}

\noindent {\bf Summary:} In general, $\Omega_S(D) = \Omega_W(D)$ and $\Omega_I
= \widehat{X_0} = \Omega_{AG}$\,.
 
\section{Cycle spaces of lower-dimensional $G_0$-orbits}
\setcounter{equation}{0}

Let us recall the setting of \cite {GM}.  For $Z=G/Q$, 
$\gamma \in Orb_Z(G_0)$ and $\kappa \in Orb_Z(K)$ its
dual, let $G(\gamma )$ be the connected component of the
identity of 
$\{ g\in G:g(\kappa )\cap \gamma \ \text{is non-empty and compact}\ \}$.  
Note that $G(\gamma )$ is
an open $K$-invariant subset of $G$ which contains the identity.
Define ${\cal C}(\gamma ):=G(\gamma )/K$.  Finally, define
${\cal C}$ as the intersection of all such cycle spaces 
as $\gamma $ ranges over $Orb_Z(G_0)$ and $Q$ ranges over all 
parabolic subgroups of $G$.

\begin {theorem}
${\cal C}=\Omega _{AG}$\,.
\end {theorem}

This result was checked in \cite {GM} for classical and hermitian exceptional
groups by means of case by case computations, and the authors of \cite{GM}
conjectured it in general.  As will be shown here, it is a
consequence of the statement 
\begin {equation*}
\Omega _W(D)=\Omega _S(D) \text{ when } D \text{ is an open }
G_0\text{--orbit in } G/B,
\end {equation*}
and of the following general result \cite[Proposition 8.1]{GM}.

\begin {proposition}
$\left ( \bigcap_{D \subset G/B \text{ open}}\ \Omega_W(D) \right ) \subset \cC$.
\end {proposition}

{\noindent {\bf Proof of Theorem.}
The 
polar $\widehat {X_0}$ in $Z=G/B$ coincides with the 
cycle space ${\cal C}_Z(\gamma _0)$, where $\gamma _0$ is the unique
closed $G_0$-orbit in $Z$. As was shown above, this agrees
with the Iwasawa domain $\Omega _I$\, which in turn is contained
in every Schubert domain $\Omega _S(D)$. Thus, for every open
$G_0$-orbit $D_0$ in $Z=G/B$ we have the inclusions

\begin {equation*}
\left ( {\bigcap}_{D \subset G/B \text{ open}}\ \Omega_W(D) \right )
\subset {\cal C}\subset {\cal C}_Z(\gamma _0)=
\widehat{X_0} =\Omega _I\subset \Omega _S(D_0)=\Omega _W(D_0).
\end {equation*}
Intersecting over all open $G_0$-orbits $D$ in $G/B$, the equalities 
\begin {equation*}
\left ( {\bigcap}_{D \subset G/B \text{ open}}\ \Omega_W(D) \right )
={\cal C}=\Omega _I=
\left ( {\bigcap}_{D \subset G/B \text{ open}}\ \Omega_W(D) \right )
\end {equation*}
are forced, and ${\cal C}=\Omega _{AG}$ is a consequence of 
$\Omega _I=\Omega _{AG}$. \hfill $\square $

\medskip
As noted in our introductory remarks, using in particular the
results of the present paper, it was shown in \cite {FH} that 
$\Omega _W(D) =\Omega _{AG}$ with the obvious exceptions in the
well-understood hermitian cases.  This is an essentially stronger
result than the above theorem on intersections. On the other
hand, it required a good deal of additional work and therefore
it is perhaps of interest that the intersection result follows
as above in a direct way from $\Omega_W(D)=\Omega _S(D)$. So,
for example, in any particular case where this latter point
was verified, the intersection theorem would be immediate (see
e.g. \cite {HW} for the case of $SL(n,\mathbb H)$).

\section{Groups of hermitian type} \label {Hermitian}
\setcounter{equation}{0}

Let $G_0$ be of hermitian type.  Write $\cB$ for the bounded symmetric
domain $G_0/K_0$ with a fixed choice of invariant complex structure.
Drop the colocation convention leading to (\ref{colocation}), so that
now the cycle space $\Omega_W(D)$ really consists of 
cycles as in \cite{W1} and \cite{WZ1}.
It has been conjectured (see \cite{WZ1}) that, whenever $D$ is an open 
$G_0$--orbit in
a complex flag manifold $Z = G/Q$, there are just two possibilities:
\begin{enumerate}
\item A certain double fibration (see \cite{WZ1})
is holomorphic, and $\Omega_W(D)$ is 
biholomorphic either to $\cB$ or to $\overline{\cB}$, or 
\item both $\Omega_W(D)$ and $\cB \times \overline{\cB}$ 
have natural biholomorphic embeddings into $G/K$, and there 
$\Omega_W(D) = \cB \times \overline{\cB}$\,.
\end{enumerate}
The first case is known (\cite{W1}, \cite{WZ1}), and the second case
has already been checked \cite{WZ1} in the cases where $G_0$ is a
classical group.

The inclusion $\Omega_W(D) \subset \cB \times \overline{\cB}$ was
proved in general (\cite{WZ1}; or see \cite{WZ2}).  
It is also known \cite{BHH} that $\cB \times \overline{\cB} = \Omega_{AG}$\,.
Combine this with $\Omega_{AG} \subset \Omega_I$ (\cite{H}; or see 
Theorem \ref{i-ag}), with 
$\Omega_W(D) = \Omega_S(D)$ (Corollary \ref{wolf=schubert}), and with 
$\Omega_I \subset \Omega_S(D)$ (compare definitions (\ref{def_schubert_domain})
and (\ref{def_iwasawa_domain})) to see that
\begin{equation}
\Omega_S(D) = \Omega_W(D) \subset (\cB \times \overline{\cB})
        = \Omega_{AG} \subset \Omega_I \subset \Omega_S(D).
\end{equation}
Now we have proved the following result.  (Also see \cite{WZ2}.)
 
\begin{theorem}\label{hermitian-case}
Let $G_0$ be a simple noncompact group of hermitian type.  Then either
{\rm (1)} a certain double fibration {\rm (see \cite{WZ1})}
is holomorphic, and $\Omega_W(D)$ is
biholomorphic to $\cB$ or to $\overline{\cB}$, or {\rm (2)}
$\Omega = G/K$ and $\Omega_W(D) = \Omega_S(D) = \Omega_I = \Omega_{AG}
= (\cB \times \overline{\cB})$.
\end{theorem}

\begin{remark} {\em
Since the above argument already uses the inclusion $\Omega_W(D) \subset
(\cB \times \overline{\cB})$ of \cite{WZ1}, it should be noted that
the construction for the proof of Corollary \ref{wolf=schubert}
can be replaced by the following treatment.  Using Corollary
\ref{special-slice}, given $p \in \text{ bd}(\Omega_W(D))$, one has an
Iwasawa--Borel subgroup $B \subset G$ and a $B$--invariant incidence
variety $A_Y$ such that $p \in A_Y \subset \Omega \setminus \Omega_W(D)$.
Since the open $B$--orbit in $\Omega$ is affine, it follows that 
$A_Y$ is contained in a $B$--invariant hypersurface $H$.  But
$\Omega_W(D) \subset \Omega_{AG} \subset \Omega_I$\,, and
$H \subset \Omega \setminus \Omega_I$ by definition of the latter.
Thus $\Omega_W(D) = \Omega_I = \Omega_{AG} = (\cB \times \overline{\cB})$.
}
\end{remark}

\section{Appendix: Lifting the cycle space to $G/K$} \label{lifting}
\setcounter{equation}{0}

As mentioned in connection with the definitions (\ref{def_cyclespace}) and
(\ref{def_liftedcyclespace}), we can view the cycle space $\Omega_W(D)$
inside $G/K$ because of

\begin{theorem}\label{same}  The projection $\pi : G/K \to G/\widetilde{K}$
restricts to a $G_0$--equivariant holomorphic cover 
$\pi: \widetilde{\Omega_W(D)} \to \Omega_W(D)$, and
$\pi: \widetilde{\Omega_W(D)} \to \Omega_W(D)$ is one to one.
\end{theorem}

We show that $\widetilde{\Omega_W(D)}$ is homeomorphic to a cell,
and then we apply \cite[Corollary 5.3]{F}.
\medskip

Without loss of generality we may assume that $G$ is simply connected.
Let $G_u$ denote the $\theta$--stable compact real form of $G$ such
that $G_u \cap G_0 = K_0$\,, connected.  Then $G_u$ is simply connected because
it is a maximal compact subgroup of the simply connected group $G$. It 
follows that 
$G_u/K_0$ is simply connected.  We view $G_u/K_0$ as a riemannian symmetric
space $M_u$\,, using the negative of the Killing form of $G_u$
for metric and $\theta|_{G_u}$ for the symmetry at $1\cdot K_0$\,.
It is connected and simply connected.

\begin{definition} \label{half}
{\rm Let $x_0$ denote the base point $1\cdot K_0 \in
G_u/K_0 = M_u$\,.  Let $L \subset T_{x_0}(M_u)$
denote the conjugate locus at $x_0$\,, all tangent vectors $\xi$ at $x_0$ such
that $d\exp_{x_0}$ is nonsingular at $t\xi$ for $0 \leqq t < 1$ but
singular at $\xi$.  Then we define
$$
\tfrac{1}{2}M_u :=
\left \{\exp_{x_0}(t\xi) \left |  \xi \in L \text{ and } 0 \leqq t
< \tfrac{1}{2} \right . \right \}.
$$
The conjugate locus $L$ and the cut locus are the same for $M_u$
\cite{C}, so $\frac{1}{2}M_u$ consists of the points in
$M_u$ at a distance from $x_0$ less than half way to
$\exp_{x_0}(L)$.  For example, if $G_0/K_0$ is a bounded symmetric
domain $\cB$, then a glance at the polysphere that sweeps out 
$M_u$ under the action of $K_0$ shows that $\frac{1}{2}M_u = \cB$.
}
\end{definition}

\begin{proposition} \label{cycle-from-half}  The lifted cycle space
$\widetilde{\Omega_W(D)} = G_0\cdot \frac{1}{2}M_u \subset G/K$.
It is $G_0$--equivariantly diffeomorphic to $(G_0/K_0) \times M_u$\,.
In particular it is homeomorphic to a cell.
\end{proposition}

\noindent {\bf Proof.}  According to \cite[Theorem 5.2.6]{FH},
the lifted cycle space $\widetilde{\Omega_W(D)} \subset G/K$ coincides with
the Akhiezer--Gindikin domain $\Omega_{AG}$.  The restricted root description
\cite{AG} of $\Omega_{AG}$ is (in our notation)
$$
\Omega_{AG} = G_0\cdot \exp(\{\xi \in \mathfrak a_u \mid |\alpha(\xi)| <
\pi/2 \ \ \forall \alpha \in \Delta(\mathfrak g, \mathfrak a)\})K/K,
$$
where $\mathfrak a_u$ is a maximal abelian subspace of
$\{\xi \in \mathfrak g_u \mid \theta(\xi) = -\xi\}$ and
$\Delta(\mathfrak g, \mathfrak a)$ is the resulting family of restricted
roots.  A glance at Definition \ref{half} shows that
$$
\tfrac{1}{2}M_u = K_0\cdot \exp(\{\xi \in \mathfrak a_u \mid
|\alpha(\xi)| < \pi/2 \ \
\forall \alpha \in \Delta(\mathfrak g, \mathfrak a)\})K/K.
$$
Thus $\widetilde{\Omega_W(D)} = \Omega_{AG} = G_0 \cdot
\exp(\{\xi \in \mathfrak a_u \mid |\alpha(\xi)| < \pi/2 \ \
\forall \alpha \in \Delta(\mathfrak g, \mathfrak a)\})K/K
= G_0\cdot \frac{1}{2}M_u$\,. That is the first assertion.

For the second assertion note that $\widetilde{\Omega_W(D)}$ fibers
$G_0$--equivariantly over $G_0/K_0$ by $gx \mapsto gK_0$ for $g \in G_0$
and $x \in \frac{1}{2}M_u$\,.  For the third assertion note that
$G_0/K_0$ and $\frac{1}{2}M_u$ are homeomorphic to cells.
\hfill $\square$
\medskip

\noindent {\bf Proof of Theorem.}  As $\widetilde{\Omega_W(D)}$ is a cell,
$H^q(\widetilde{\Omega_W(D)};\Z) = 0$ for $q > 0$ and the Euler characteristic
$\chi(\widetilde{\Omega_W(D)}) = 1$.  Let $\gamma$ be a covering transformation
for $\pi : \widetilde{\Omega_W(D)} \to \Omega_W(D)$ and let $n$ be its order.
Then the cyclic group $\langle \gamma \rangle$ acts freely on
$\widetilde{\Omega_W(D)}$ and the quotient manifold
$\widetilde{\Omega_W(D)}/\langle \gamma \rangle$ has Euler characteristic
$\chi(\widetilde{\Omega_W(D)})/n$ \cite[Corollary 5.3]{F}, so $n = 1$.
Now the covering group of $\pi : \widetilde{\Omega_W(D)} \to \Omega_W(D)$
is trivial, so $\pi$ is one to one.  \hfill $\square$

\centerline{\begin{tabular}{ll}
ATH: & JAW: \\
Fakult\" at f\" ur Mathematik & Department of Mathematics \\
Ruhr--Universit\" at Bochum & University of California \\
D-44780 Bochum, Germany & Berkeley, California 94720--3840, U.S.A. \\
                               &                                      \\
{\tt ahuck@cplx.ruhr-uni-bochum.de} & {\tt jawolf@math.berkeley.edu}
\end{tabular}}

\enddocument
\end